\documentclass[11pt]{article}
\usepackage{amsmath,amsfonts,amssymb,amsthm,mathrsfs}
\usepackage{wasysym}
\usepackage{minitoc}
\usepackage{endnotes}
\usepackage{manfnt}
\usepackage[dvipsnames]{xcolor}
\usepackage[a4paper,vmargin={1cm,1cm},hmargin={2cm,2cm}]{geometry}

\usepackage{graphicx,graphics}
\usepackage{epsfig}
\usepackage{latexsym}
\usepackage[applemac]{inputenc}
\usepackage{ae,aecompl}

\usepackage[english]{babel}
 \usepackage[colorlinks=true]{hyperref}
\usepackage{pstricks}
\usepackage{enumerate}
\renewcommand{\leq}{\leqslant}
\renewcommand{\geq}{\geqslant}
\usepackage{skull}
\usepackage{mathpazo}
\usepackage[font=sf, labelfont={sf,bf}, margin=0.5cm]{caption}

\usepackage{titlesec}
\usepackage{titletoc}
\usepackage{color}

\titleformat{\part}[display]{\beaupetit}{\Huge\textsf{Part \textsc{\Roman{part}} : }    }{0pt}{}[]

\DeclareFixedFont{\chapterNumberFont}{OT1}{fsk
}{b}{n}{5cm} 
\DeclareFixedFont{\classik}{OT1}{ptm}{b}{sc}{1cm} 
\DeclareFixedFont{\chapterTitleFont}{T1}{ftp}{b}{n}{1cm} 
\definecolor{gris75}{gray}{0.75} 
\definecolor{gris25}{gray}{0.15}

\titleformat{\section}[block]{\Large\sffamily}{\noindent \bfseries \thesection}{1em}{} 
\titleformat{\subsection}[block]{\large\sffamily}{\thesubsection}{1em}{} 
\titleformat{\subsubsection}[block]{\sffamily}{}{}{} 
\titleformat{\paragraph}[runin]{\sffamily}{}{}{} 
\titlespacing{\section} {0pc}{3.5ex plus .1ex minus .2ex}{1.5ex minus .1ex}

\titleformat{\chapter}[hang]{\bfseries\Large\chapterTitleFont}{\textsf{Chapter \textsc{\Roman{chapter}} : }  }{0pt}{}[\titlerule 
]

  \DeclareFixedFont{\beaupetit}{T1}{ftp}{b}{n}{2cm} 

\usepackage{thmtools}
\declaretheorem[thmbox=M,name=Theorem]{theorem}
\declaretheorem[name=Lemma,sibling=theorem]{lemma}

\def\build#1_#2^#3{\mathrel{
\mathop{\kern 0pt#1}\limits_{#2}^{#3}}}

             \title{Yet another proof of the strong law of large numbers}
             \author{Nicolas Curien}
             \date{}
             
\begin{document}
             
             \maketitle
             \vspace{-1cm}\begin{abstract} We give a short proof of the strong law of large numbers based on duality for random walk.\end{abstract}
\begin{figure}[!h]
 \begin{center}
 \includegraphics[height=3.5cm]{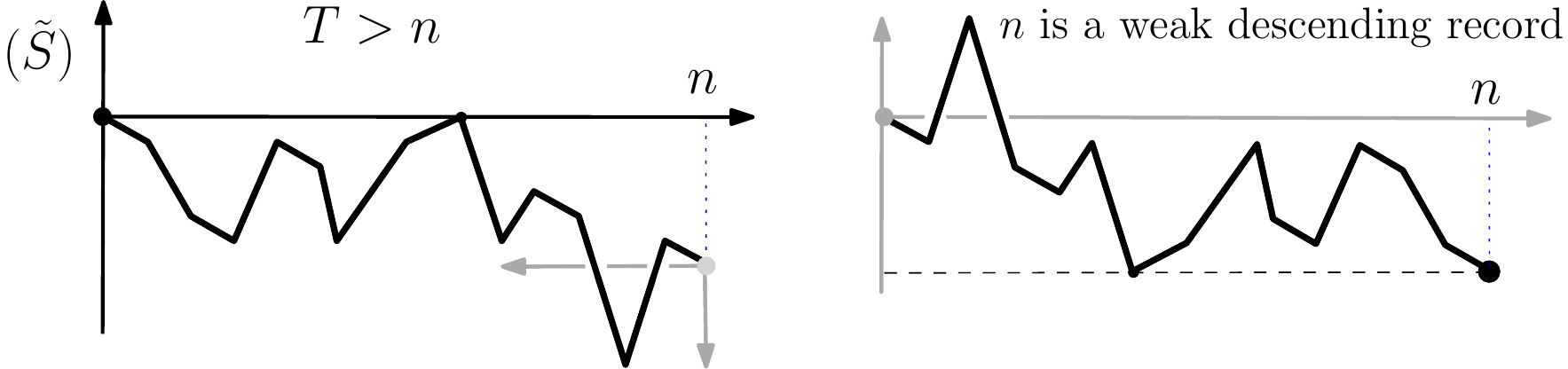}
 \caption{ \label{fig:duality} Time and space reversal shows that $\mathbb{P}(T>n) = \mathbb{P}(n \mbox{ is a descending record time}).$}
 \end{center}
 \end{figure}

             \thispagestyle{empty}
 Let $X_{1}, X_{2}, \dots $ be i.i.d.~random variables with finite expectation and let $S_{n} = X_{1}+ \dots + X_{n}$ for $n \geq 0$ be the corresponding random walk. Kolmogorov's strong law of large numbers says that $ n^{-1}S_{n} \to  \mathbb{E}[X]$ almost surely as $n \to \infty$. Clearly, it is a consequence of the lemma:
             
             \begin{lemma} Let $X_{1}, X_{2}, \dots $ be i.i.d.~r.v.~with $ \mathbb{E}[X]>0$. Then $ \inf_{n \geq 0} (X_{1} + \cdots + X_{n})$ is finite a.s.
             \end{lemma}
\noindent \textbf{Proof.}
 \textsc{Step 1. Bounding the increments from above.} Choose $C>0$ large enough so that by dominated convergence $ \mathbb{E}[X \mathbf{1}_{X<C}] > 0$. We will show that the random walk $\tilde{S}_{n} = X_{1} \mathbf{1}_{X_{1}<C} + \dots + X_{n} \mathbf{1}_{ X_{n}<C}$ is a.s.~bounded from below which is sufficient to prove the lemma.\\
\textsc{Step 2. Duality.} For every $n \geq 0$ we have the equality in law       $ (0=\tilde{S}_{0},\tilde S_{1}, \dots , \tilde S_{n}) {=} (\tilde S_{n}- \tilde S_{n}, \tilde S_{n}- \tilde S_{n-1}, \dots , \tilde S_{n}-\tilde S_{1}, \tilde S_{n}-\tilde S_0).$
Let $T = \inf\{i \geq 0 : \tilde{S}_{i}>0\}$ be the first hitting time of the positive axis by the walk and recall that a time $n\geq 0$ is a weak descending record time if and only if $ \tilde{S}_{n} = \min_{0 \leq k \leq n} \tilde S_{k}.$ By applying the above equality in law we deduce (see Figure \ref{fig:duality}) that 
$$ \mbox{ for all }n \geq 0, \quad  \mathbb{P}(T>n) = \mathbb{P}(n \mbox{ is a weak descending record time}).$$
Summing over $n\geq 0$, we get that $ \mathbb{E}[T] = \mathbb{E}[\# \mbox{ weak descending record times}]$ and the proof is complete if we prove that  $ \mathbb{E}[T]< \infty$ since this implies that almost surely there is a finite number of weak descending records for $ \tilde{S}$, hence the walk is bounded from below a.s.\\
\textsc{Step 3. Optional sampling theorem.} To prove  $ \mathbb{E}[T]<\infty$, consider the standard martingale $$ M_{n}~=~\tilde{S}_{n} - \mathbb{E}[X \mathbf{1}_{X<C}]n, \quad \mbox{ for }n \geq0$$ (for the filtration generated by the $X_{i}$'s) and apply the optional sampling theorem to the stopping time $n \wedge T$ to deduce that 
$$ 0=\mathbb{E}[M_{n\wedge T}] \quad \mbox{ or in other words } \quad \mathbb{E}[X \mathbf{1}_{X<C}] \mathbb{E}[n \wedge T] = \mathbb{E}[\tilde{S}_{n \wedge T}].$$
Since the increments of $\tilde{S}$ are bounded above by $C$, the right-hand side of the last display is bounded by $C$ as well. Letting $n \to \infty$, by monotone convergence we deduce that the expectation of $T$ is finite. Et voilà. \qed \medskip 

\noindent \textbf{Acknowledgments:} We thank Yuval Peres and the probability team at Orsay for helpful feedback.
\end{document}